\documentclass[12pt]{amsart}

\usepackage{amsfonts, amssymb, amsgen, amsbsy, amstext, amsopn, amsmath,
latexsym, indentfirst,color,longtable,verbatim}
\usepackage[utf8]{inputenc}
\usepackage[T1]{fontenc}
\usepackage[english]{babel}
\usepackage{verbatim} 

\renewcommand{\H}{\mathbb{H}}
\newcommand{\B}{\mathbb{B}}

\newcommand{\G}{\mathbb{G}}

\newcommand{\N}{\mathbb{N}}

\newcommand{\R}{\mathbb{R}}

\newcommand{\cC}{\mathcal{C}}

\newcommand{\cF}{\mathcal{F}}

\newcommand{\cH}{\mathcal{H}}

\newcommand{\cP}{\mathcal{P}}

\newcommand{\cS}{\mathcal{S}}

\newcommand{\ep}{\varepsilon}

\newcommand{\ph}{\varphi}

\newcommand{\bcup}{\bigcup}

\newcommand{\lan}{\langle}
\newcommand{\ran}{\rangle}

\newcommand{\der}{\partial}

\newcommand{\res}{\mbox{\LARGE{$\llcorner$}}}
\newcommand{\Rd}{{\mathcal F}_HE}
\newcommand{\Per}{|\partial_HE|}

\newcommand{\LLb}{{\bigg\lbrace}}
\newcommand{\RRb}{{\bigg\rbrace}}

\newcommand{\ls}{\mbox{\large $($}}
\newcommand{\rs}{\mbox{\large $)$}}
\newcommand{\lls}{\big(}
\newcommand{\rrs}{\big)}

\newcommand{\diam}{\mbox{\rm diam}}

\newcommand{\ds}{\displaystyle}
\newcommand{\beqas}{\begin{eqnarray*}}
\newcommand{\eeqas}{\end{eqnarray*}}
\newcommand{\beqa}{\begin{eqnarray}}
\newcommand{\eeqa}{\end{eqnarray}}
\newcommand{\beq}{\begin{equation}}
\newcommand{\eeq}{\end{equation}}
\newcommand{\bce}{\begin{center}}
\newcommand{\ece}{\end{center}}

\newcommand{\pa}[1]{\left( #1 \right)}               
\newcommand{\set}[1]{\left\{ #1 \right\}}            
\newcommand{\qandq}{\quad\mbox{and}\quad}
\newcommand{\qforeveryq}{\quad\mbox{for every}\quad}

\newtheorem{The}{Theorem}[section]
\newtheorem{Lem}{Lemma}[section]
\newtheorem{Def}{Definition}[section]
\newtheorem{Rem}{Remark}[section]
\newtheorem{Pro}[The]{Proposition}
\newtheorem{Cor}{Corollary}
\newtheorem{Exa}[The]{Example}

\newcommand{\bt}{\begin{The}}
\newcommand{\et}{\end{The}}
\newcommand{\bl}{\begin{Lem}}
\newcommand{\el}{\end{Lem}}
\newcommand{\bd}{\begin{Def}}
\newcommand{\ed}{\end{Def}}
\newcommand{\br}{\begin{Rem}}
\newcommand{\er}{\end{Rem}}
\newcommand{\bpr}{\begin{Pro}}
\newcommand{\epr}{\end{Pro}}
\newcommand{\bc}{\begin{Cor}}
\newcommand{\ec}{\end{Cor}}
\newcommand{\bex}{\begin{Exa}}
\newcommand{\eex}{\end{Exa}}

\newcommand{\bR}{\ensuremath{\mathbb{R}}}

\begin{document}

\title
[Some remarks on densities in the Heisenberg group]
{Some remarks on densities in the Heisenberg group}
\author{Valentino Magnani}
\address{Valentino Magnani,  Dipartimento di Matematica \\
Largo Pontecorvo 5 \\ I-56127, Pisa}
\email{magnani@dm.unipi.it}
%
%
%
%
\date{\today}

\begin{abstract}
We observe that upper densities and spherical Federer densities may differ on all two dimensional surfaces of the sub-Riemannian Heisenberg group. 
This provides an entire class of intrinsic rectifiable sets having upper density strictly less than one.
\end{abstract}

\maketitle

\section{Introduction}

Area formulae in homogeneous groups can be seen in the wide project to develop Geometric Measure Theory in a non-Euclidean framework. Here a basic difficulty is the possible difference between Hausdorff dimension and topological dimension of smooth submanifolds. It turns out that such submanifolds cannot rectifiable, according to the metric notion of rectifiability, \cite{Federer69}.
This problem already appears in the Heisenberg group, that is the simplest model of noncommutative homogeneous group, \cite{Stein93}.

Measure theoretic area formulae are the key tool to overcome the lack of Lipschitz parametrizations. 
If $\cS^\alpha$ is the $\alpha$-dimensional spherical measure, then
\begin{equation}\label{eq:spharea}
 \mu(B)=\int_B \theta^\alpha(\mu,x)\,d\cS^\alpha(x)
\end{equation}
for a Borel regular measure $\mu$ on a metric space $X$, a Borel set $B\subset X$ and where clearly $\mu\res B$ is absolutely continuous
with respect to $\cS^\alpha\res B$, see \cite{Mag30} for more information. 
The point of \eqref{eq:spharea} is the explicit formula of the {\em spherical Federer density} $\theta^\alpha(\mu,\cdot)$, namely, 
\beq
\theta^\alpha(\mu,x)=\inf_{\ep>0}\,\sup\LLb\frac{\mu(\B)}{c_\alpha \diam(\B)^\alpha}: x\in \B, \B \ \mbox{is a closed ball},\ \diam(\B)<\ep\RRb.
\eeq
Thus, densities of measures are strictly related to measure theoretic area formulae. 
On the other hand, $\theta^\alpha(\mu,x)$ is a kind of ``non-centered upper density'', that may differ from the standard {\em upper $\alpha$-density}
\[
\Theta^{*\alpha}(\mu,x)=\limsup_{r\to0^+}\frac{\mu\lls\B(x,r)\rrs}{c_\alpha r^\alpha}.
\]
We refer to 2.10.19 of \cite{Federer69} for more information on upper and lower $\alpha$-densities in metric spaces. When they coincide, 
their common value is the $\alpha$-density, denoted by $\Theta^\alpha(\mu,\cdot)$.
%
\begin{comment}
General Federer densities with respect to a family of sets and a size function have been recently introduced in \cite{Mag30}.
\end{comment}
%

In the important case of rectifiable measures, where $E\subset\R^n$ is $k$-rectifiable, $\cH^k_{|\cdot|}$ is the Euclidean Hausdorff
measure and setting $\mu=\cH^k_{|\cdot|}\res E$, we have  
\beq\label{eq:Hktheta}
\Theta^k(\cH^k_{|\cdot|}\res E,\cdot)=\theta^k(\cH^k_{|\cdot|}\res E,\cdot)=1
\eeq
%
\begin{comment}
...that here it coincides with the Euclidean spherical measure $\cS^k_{|\cdot|}$
In fact, the spherical Federer density $\theta^k(\cH^k_{|\cdot|}\res E,\cdot)$ equals one a.e., as it immediately follows from \eqref{eq:spharea}.
\end{comment}
%
$\cH^k_{|\cdot|}$ a.e.\ in $E$, see for instance 3.2.19 of \cite{Federer69}. This implication holds for more general rectifiable sets in metric spaces,
\cite{Kir94}. Conversely, in Euclidean space it is well known that the validity of \eqref{eq:Hktheta} a.e.\ implies the $m$-rectifiability of $E$,  \cite{Mattila75}. 
Moreover, the only existence of the $m$-density a.e.\ implies the same $m$-rectifiability, \cite{Preiss87}.
A recent account on densities and rectifiability in the Euclidean space can be  found in \cite{DeLellis08}.
%
\begin{comment} 
We simply observe that $\cH^k=\cS^k$ and the application of \eqref{eq:spharea} to $\mu=\cH^k\res E$ immediately gives  $\theta^k(\cH^1\res E,x)=1$ for $\cH^k$ a.e.\ $x\in E$. 
 However in Heisenberg groups and more general stratified groups rectifiable sets in the classical sense may not have nontrivial examples ...ambkir mag

The previous equality can be extended to rectifiable sets of metric spaces, \cite{Kir94}.

\end{comment}
%

Natural notions of ``intrinsic rectifiability'' in Heisenberg groups and general stratified groups have been also studied, proving an analogous role of the classical rectifiability, \cite{FSSC3}, \cite{FSSC4}, \cite{Mag5}. Thus, we may ask at which extent we can expect to find deep relationships between rectifiability and densities in these groups.

A first promising result is by Mattila, Serapioni and Serra Cassano, who characterize the intrinsic rectifiability in Heisenberg groups by the a.e.\ existence of suitable tangent subgroups of fixed dimension, \cite{MSSC10}.
More recently, the Marstrand's density theorem in Heisenberg groups have been proved by Chousionis and Tyson, \cite{ChoTys}. 
 
Our question is the following: is it reasonable to expect \eqref{eq:Hktheta} to persist also when $E$ is an intrinsic rectifiable set of the first Heisenberg group?
These sets, called $\H$-rectifiable sets, are essentially countable unions of level sets of functions having nonvanishing differential along the directions of the horizontal subbundle of the Heisenberg group, \cite{FSSC3}. Notice that the defining functions are not necessarily differentiable in the classical 
sense, \cite{Mag6}. In fact, $\H$-rectifiable sets form a strictly larger class than smooth surfaces,
 \cite{Balogh03}, \cite{FSSC3}, \cite{KirSer04}.
 
The following theorem answers the above question, when the sub-Riemannian distance $\rho$, in short SR distance, is fixed in the
first Heisenberg group $\H$.
\bt\label{th:densless1}
There exist a geometric constant $0<\gamma_\rho<1$, only depending on $\rho$, such that for each $\H$-rectifiable set $\Sigma\subset \H$ we have 
\beq\label{eq:ineqS}
\Theta^{*3}(\cH_\rho^3\res\Sigma,x)\le\gamma_\rho
\eeq
for $\cH_\rho^3$-a.e.\ $x\in\Sigma$, where $\cH_\rho^3$ is the Hausdorff measure constructed by $\rho$.
\et
This theorem is strictly related to the shape of the metric unit ball $\B_\rho$, see \eqref{eq:Brho}. In fact, the geometric constant $\gamma_\rho$ is a quotient between the areas of two suitable slices of $\B_\rho$, according to \eqref{eq:betarho}, \eqref{eq:beta0rho} and \eqref{eq:strict1}. 
In view of the previous comments, Theorem~\ref{th:densless1} shows that the SR distance is not suitable to develop analogues of the classical rectifiability theorems.
%
\begin{comment}
It is true that the difference between these densities was already shown only for nonhorizontal curves of the Heisenberg group with respect to the sub-Riemannian distance, \cite{Mag30}, but it is not yet clear whether these sets can be considered in some sense intrinsically rectifiable. Certainly they are not rectifiable in the sense of Rumin.

--------

In the case the sub-Riemannian distance $\rho$ is replaced by a distance $d$, whose unit ball is convex,
then under the assumptions of Theorem~\ref{th:densless1} we would obtain $\Theta^3(\cS_d^3\res\Sigma,x)=1$ for $\cS_d^3$-a.e.\ $x\in\Sigma$, from the results of \cite{Mag31}. 
\end{comment}
%
However, a few intriguing questions are now in order. 
If we choose any other homogeneous distance of $\H$, see Definition~\ref{def:homdim}, what can we say about the existence 
of the density $\Theta^{3}(\cH_d^3\res\Sigma,\cdot)$ and its value? It is also unclear whether the upper density $\Theta^{*3}(\cH_d^3\res\Sigma,\cdot)$ is either one or strictly less than one $\cH_d^3$-a.e.\ in $\Sigma$.

The proof of Theorem~\ref{th:densless1} follows from the fact that the 3-density for the spherical measure $\cS^3_\rho\res\Sigma$ actually exists, that is
\beq\label{eq:thetaS3}
\Theta^{3}(\cS_\rho^3\res\Sigma,x)=\gamma_\rho
\eeq
for $\cS^3_\rho$-a.e.\ $x\in\Sigma$. This follows from the a.e.\ blow-up of perimeter measure and its integral representation with respect to $\cS^3_\rho$, see  \cite{Mag31}. 

We wish to stress that the strict inequality $\gamma_\rho<1$ is possible due to the nonconvex shape of metric ball $\B_\rho$. This special SR phenomenon could not appear in $\R^n$ equipped with the usual commutative operations, 
since here any homogeneous distance gives a Banach norm and all metric balls of a Banach space are obviously convex sets.
In fact, in any finite dimensional Banach space Theorem~\ref{th:densless1} cannot hold,
as a consequence of the general results of \cite{Kir94}.

The equality \eqref{eq:thetaS3} also affects the relationship between the spherical measure $\cS^3_\rho$ and 
the centered Hausdorff measure $\cC^3_\rho$. The latter, especially known in Fractal Geometry, was introduced by Saint Raymond and Tricot, \cite{RayTri88}.
It is a variant of the spherical measure, that could be seen as ``dual'' of the packing measure, in view of its relationship to upper density as the packing measure has with its lower density, see Theorem~1.1 of \cite{RayTri88} and \cite{Edg2007} for more information. We have the following theorem.
\bt\label{th:SlessC3}
For each $\H$-rectifiable set $\Sigma\subset\H$ with $\cS^3$ positive measure there holds
\beq\label{eq:S3C3}
\cS^3_\rho(\Sigma)<\cC^3_\rho(\Sigma).
\eeq
In particular, we have $\cS^3_\rho<\cC^3_\rho$ as measures.
\et
The proof of this theorem also relies on the results of \cite{Mag31}, joined with a recent measure theoretic area formula, proved by Franchi, Serapioni and Serra Cassano \cite{FSSC8}. It is worth to compare Theorem~\ref{th:SlessC3} with Euclidean results, where spherical measures and centered Hausdorff measures always coincide on rectifiable sets, \cite{RayTri88}. 
Even in $\H$, if we change the dimension into that of the Haar measure, we get $\cS^4_\rho=\cC^4_\rho$ and this holds for more general groups, \cite{FSSC8}. Moreover, in the assumptions of Theorem~\ref{th:SlessC3}, replacing $\rho$ by any homogeneous distance $d$ whose unit ball is convex, we obtain
the equality $\cC_d^3(\Sigma)=\cS_d^3(\Sigma)$, see Remark~\ref{re:cCcS3}. 
These facts show how Theorem~\ref{th:SlessC3} provides another unexpected feature of intrinsic rectifiable sets of $\H$, when seen through the SR distance.

Some additional questions are still to be understood. In fact, the same arguments for the proof of Theorem~\ref{th:densless1}, joined with Theorem~5.2 of \cite{Mag31},
would also imply that
\beq\label{eq:thetaS3d}
\Theta^{3}(\cS_d^3\res\Sigma,x)=1
\eeq
for $\cS^3_d$-a.e.\ $x\in\Sigma$, where $\Sigma$ is $\H$-rectifiable and the unit ball with respect to $d$ is convex.
Since we may choose an $\H$-rectifiable set $\Sigma_0$ that is not rectifiable, \cite{KirSer04}, it turns out that the metric space $(\H,d)$ possesses an unrectifiable set $\Sigma_0$ whose 3-density with respect to $\cS^3_d\res\Sigma$ is equal to one a.e.\ in $\Sigma$.

To the author's knowledge, the question of finding a metric space with a subset of $k$-density equal to one a.e.\ and that is not rectifiable
is still an important open question when $k>1$. For $k=1$, this problem is well known, \cite{PreissTiser92}.
Here the $k$-density refers to the $k$-dimensional Hausdorff measure. 
For this reason, the previous example does not answer the open question, being the exact formula between $\cS^3_d\res\Sigma$ and $\cH^3_d\res\Sigma$ in $\H$ still unknown, even for smooth 2-dimensional submanifolds.

%
%
%
\section{Basic definitions and proofs}\label{SectHeis}

We consider the three dimensional Heisenberg group $\H$, that is represented by $\R^3$  equipped with the Lie group operation
\[
(x,y,t)\cdot (x',y',t')=(x+x',y+y',t+t'-2xy'+2x'y)
\]
for all $(x,y,t),(x',y',t')\in\H$ and the left invariant vector fields
\[
X=\der_x+2y\der_t, \quad Y=\der_y-2x\der_t\quad\mbox{and}\quad T=\der_t\,.
\]
This system of coordinates also yields an auxiliary Euclidean structure on $\H$, where $|\cdot|$ denotes the Euclidean norm of $\R^3$,
that is automatically induced on $\H$.

Now we introduce the distance that makes $\H$ an SR manifold.
It is well known that $\H$ is pathwise connected by {\em horizontal curves}, namely absolutely continuous curves
$\gamma:[0,1]\to\H$ such that $\dot\gamma(t)$ is a linear combination of $X(\gamma(t))$ and $Y(\gamma(t))$ for a.e.\ $t\in[0,1]$.
This permits us to define the following SR distance 
\[
\rho(p,q)=\inf\set{ \int_0^1 |\dot\gamma(t)|\,dt: \ \gamma\ \mbox{is horizontal and connects $p$ with $q$}}
\]
between any couple of points $p,q\in\H$, where $|\cdot|$ is the left invariant Riemannian metric that makes the vector fields
$X,Y,Z$ everywhere orthonormal. Explicit formulae for geodesics with respect to $\rho$ and the equations for the the boundary of the SR unit ball 
\beq\label{eq:Brho}
\B_\rho=\{p\in \H: \rho(p,0)\le1\}
\eeq
are well known facts, see for instance \cite{Monti00} and \cite{Bell96}. Precisely, the boundary of $\B_\rho$ is the image of the mapping 
\[
\Phi\colon[0,2\pi]\times[-2\pi,2\pi]\to\bR^3,
\]
whose components are defined as follows
\beqa\label{eq:unitrho}
\left\{\begin{array}{l}
\Phi_1(\theta,\varphi)=\ds\frac{\cos\theta(1-\cos\varphi)+\sin\theta\sin\varphi}\varphi\\
\Phi_2(\theta,\varphi)=\ds\frac{-\sin\theta(1-\cos\varphi)+\cos\theta\sin\varphi}\varphi\\
\Phi_3(\theta,\varphi)=\ds2\left(\frac{\varphi-\sin\varphi}{\varphi^2}\right)
\end{array}\right.\,,
\eeqa
where we understand $\Phi_1(\theta,0)=\sin\theta$, $\Phi_2(\theta,0)=\cos\theta$ and $\Phi_3(\theta,0)=0$.
%
\begin{comment}
{\scriptsize L'uguaglianza
\[
\inf\set{ \int_0^1 |\dot\gamma(t)|\,dt: \ \gamma\ \mbox{horiz. and it connects $x$ with $y$}}=
\inf\set{ \sqrt{\int_0^1 |\dot\gamma(t)|^2\,dt}: \ \gamma\ \mbox{horiz. and it connects $x$ with $y$}}
\]
\`e provata nel lavoro di Lanconelli, pubblicato nel Seminario Matematico di Bologna.
Per tale uguaglianza utilizza le subunit curves e la nozione di distanza come tempo minimo con vincolo sulle velocit\`a.
Tale idee possono essere gi\`a presenti nel lavoro di Jerison e Sanchez-Calle ``subelliptic, second order differential operators'', in Complex analysis III,
cerca anche nel lavoro di Fefferman e Phong.
}
\end{comment}
%
The precise shape of the unit ball $\B_\rho$ will play an important role in the subsequent computation.

%
%
%
The distance $\rho$ allows us to construct the associated spherical measure and Hausdorff measure.
These measures can be construced in a general metric space.
\begin{Def}[Hausdorff measures]\label{def:hausdorffsph}\rm 
Let $X$ be a metric space equipped with a distance $d$. Let $\cF\subset\cP(X)$ be a nonempty class of closed subsets of $X$, let $\alpha>0$ and $c_\alpha>0$. 
If $\delta>0$ and $E\subset X$, then we define
\begin{equation*}
\phi_\delta(E)=\inf \set{\sum_{j=0}^\infty c_\alpha\, \diam(B_j)^\alpha: E\subset \bcup_{j\in\N} B_j, \diam(B_j)\le\delta,\  B_j\in\cF}\,.
\end{equation*}
When $\cF$ is the family of closed balls, then the {\em $\alpha$-dimensional spherical measure} is defined by $\cS^\alpha(E)=\sup_{\delta>0}\phi_\delta(E)$,
for every $E\subset\G$.  When $\cF$ is the family of closed subsets, then the {\em $\alpha$-dimensional Hausdorff measure} is defined by $\cH^\alpha(E)=\sup_{\delta>0}\phi_\delta(E)$.
In the case $c_\alpha=2^{-\alpha}$, we will use the notation $\cS_0^\alpha$ and $\cH^\alpha_0$.
\end{Def}
We define the spherical measure  $\cS^3_\rho$ and the Hausdorff measure $\cH^3_\rho$ by setting 
\beq\label{eq:betarho}
c_3=\beta_\rho\,2^{-3}\qandq \beta_\rho=\max_{w\in \B_\rho} \cH_{|\cdot|}^2(\B_\rho\cap wN_0),
\eeq
where we have set $N_0=\set{(x,0,z)\in\R^3}$. We also introduce the geometric constant
\beq\label{eq:beta0rho}
\beta^0_\rho=\cH^2_{|\cdot|}(\B_\rho\cap N_0).
\eeq
We will also consider a more general class of distances than $\rho$, that are compatible with the algebraic structure of $\H$.
%
%
\begin{Def}[Dilations and homogeneous distances]\label{def:homdim}\rm 
For each $s>0$, a {\em dilation} $\delta_s:\H\to\H$ is defined as 
\[
\delta_s(x,y,z)=(sx,sy,s^2z)\qforeveryq (x,y,z)\in\H.
\]
A distance $d$ in $\H$ is {\em homogeneous} if it is continuous, and for each $q,w,u\in\H$ and $r>0$ we have
$d(qw,qu)=d(w,u)$ and $d(\delta_rw,\delta_ru)=rd(w,u)$.
\ed
\br\rm
It can be shown that $\rho$ is a homogeneous distance. 
\er
%
%
\begin{Def}[Centered Hausdorff measure]\label{def:chausdorff}\rm 
Let $X$ be a metric space equipped with a distance $d$. We fix $\alpha>0$ and $c_\alpha>0$
and denote by $\cF_b$ the family of closed balls in $X$.
If $\delta>0$ and $S\subset X$, then we define
\[
\cC^\alpha_\delta(S)=\inf \set{\sum_{j=0}^\infty c_\alpha\,\diam(B_j)^\alpha:S\subset \bcup_{j\in\N} B_j,\ \mbox{$B_j\in\cF_b$ is centered in $S$, $\diam(B_j)\le\delta$} }
\]
and $\tilde{\cC}^\alpha(S)=\sup_{\delta>0}\cC^\alpha_\delta(S)$. Since $\tilde\cC^\alpha$ may not be a measure, see \cite{RayTri88},
we finally define
\[
\cC^\alpha(E)=\sup\set{\tilde{\cC}^\alpha(S): S\subset E }
\]
for every $E\subset X$. According to \cite{RayTri88}, $\cC^\alpha$ is a measure, called the $\alpha$-dimensional {\em centered Hausdorff measure}.
When $c_\alpha=2^{-\alpha}$, we use the notation $\cC_0^\alpha$.
\end{Def}
For the 3-dimensional centered Hausdorff measure $\cC_d^3$ with respect to a homogeneous distance $d$ of $\H$ 
we use the same constant $c_3=\beta_\rho\,2^{-3}$ of the spherical measure $\cS^3_\rho$
%
%
%
%

\subsection{Proof of Theorem~\ref{th:densless1}}
Following notation and definitions of \cite{Mag31}, $\Sigma$ is contained in a countable union of portions of {\em reduced boundaries} of $h$-finite perimeter sets,
up to a $\cS^3_\rho$ negligible set. Thus, it is not restrictive to assume that 
\[
\Sigma=\Omega\cap\Rd,
\]
where $\Omega\subset\H$ is an open set, $E$ is an $h$-finite perimeter of $\H$ and $\Rd$ is its reduced boundary.
By Theorem~1.2 of \cite{Mag31}, we have
\begin{equation}\label{eq:perimBeta}
\Per\res\Omega=\beta(\rho,\nu_E)\, \cS_0^3\res(\Omega\cap\Rd)=\beta(\rho,\nu_E)\, \cS_0^3\res\Sigma
\end{equation}
and for each $x\in\Sigma$ there holds
\[
\beta\lls \rho,\nu_E(x)\rrs=\max_{z\in \B_\rho} \cH_{|\cdot|}^2\ls\B_\rho(z,1)\cap N\lls\nu_E(x)\rrs\rs
\]
where $N\lls\nu_E(x)\rrs$ is the vertical plane orthogonal to the horizontal normal $\nu_E(x)$. 
The rotational symmetry of $\B_\rho$ implies that $\beta(\rho,\cdot)$ is constant and the fact that translations preserve
the Euclidean Hausdorff measure $\cH^2_{|\cdot|}$ between translated vertical planes gives
\[
\beta_\rho=\beta(\rho,\cdot).
\]
For the same reasons, the function 
\[
v\to \cH^2_{|\cdot|}\lls\B_\rho\cap N(v)\rrs=\beta^0_\rho
\]
is constant, as $v$ varies in $\R^2\times\{0\}$. As a result, by \eqref{eq:perimBeta}, we get $\Per\res\Omega=\cS_\rho^3\res\Sigma$. 
At $\Per$-a.e.\ $x\in\Sigma$, by Theorem~3.1 of \cite{FSSC5} we get
\[
\lim_{r\to0^+}\frac{|\der_HE|(\B(x,r))}{r^3}=\cH^2_{|\cdot|}(\B_\rho\cap N_0)=\beta^0_\rho.
\]
It follows that
\[
\Theta^{*3}(\cS^3_\rho\res\Sigma,x)=\limsup_{r\to0^+}\frac{|\der_HE|(\B(x,r))}{\beta_\rho r^3}=\Theta^3(\cS^3_\rho\res\Sigma,x)=\frac{\beta^0_\rho}{\beta_\rho}.
\]
Finally, our claim follows by checking that $\beta^0_\rho<\beta_\rho$. This is a direct computation that can be carried out by equations \eqref{eq:unitrho}.
Since $ \Phi_1^2(\theta,\varphi)+\Phi_2^2(\theta,\varphi)=(2-2\cos\varphi)/\varphi^2$
one realizes that the intersection of $\der\B_\rho$ with $\set{(x,0,z): x,z\ge0}$ can be parametrized by the curve
\beq\label{eq:phitil}
\tilde\Phi(\ph)=\pa{\frac{\sqrt{2-2\cos\varphi}}\varphi,0, \frac{2\varphi-2\sin\varphi}{\varphi^2}}
\eeq
defined on $[0,2\pi]$, where we understand $\tilde\Phi(0)=(1,0,0)$.
It follows that the closed upper half of $\der\B_\rho$ is the image of the mapping $F\colon[0,2\pi]\times[0,2\pi]\to\bR^3$ defined as
\beqa\label{ball}
\left\{\begin{array}{l}
\ds F_1(\psi,\varphi)=\frac{\sqrt{2-2\cos\varphi}}\varphi\cos\psi \\
\ds F_2(\psi,\varphi)=\frac{\sqrt{2-2\cos\varphi}}\varphi\sin\psi \\
\ds F_3(\psi,\varphi)=2\left(\frac{\varphi-\sin\varphi}{\varphi^2}\right)
\end{array}\right.\,,
\eeqa
where $F_1(\psi,0)=\cos\psi$, $F_2(\psi,0)=\sin\psi$ and $F_3(\psi,0)=0$.
From both the rotational and the antipodal symmetry of $\B_\rho$, if we set
$A_\rho=\B_\rho\cap\set{(x,0,z)\in\R^3:x,y\ge0}$
and take into account \eqref{eq:beta0rho}, then Gauss-Green's theorem applied to the curve $\tilde\Phi$ of \eqref{eq:phitil} gives
\begin{equation}\label{eq:0sect}
\frac{\beta^0_\rho}{4}=\cH^2(A_\rho)=2\int_0^{2\pi}\frac{\sqrt{2-2\cos\varphi}}{\varphi^4}\,\ls 2\sin\varphi-\varphi\cos\varphi-\varphi\rs d\varphi.
\end{equation}
Intersecting $\B_\rho$ with the plane $N_1=\set{ (x,y,z)\in\bR^3\colon y=2\sqrt{2}/3\pi}$, we get the curve $\tilde F$ on $[0,3\pi/2]$ defined as 
%
\begin{comment} 
Fix $(\psi_0,\varphi_0)\in[0,2\pi[\times[0,2\pi[$ such that $F_2(\psi_0,\varphi_0)=2\sqrt{2}/3\pi$ and $F_1(\psi_0,\varphi_0)\geq0$, where $F$ is the mapping from~\eqref{ball}. Using, as before, the fact that
$$
F_1^2(\psi,\varphi)+F_2^2(\psi,\varphi)=\frac{2-2\cos\varphi}{\varphi^2},
$$
when $\varphi\neq0$, we observe that
$$
F_1(\psi_0,\varphi_0)=\sqrt{\frac{2-2\cos\varphi_0}{\varphi_0^2}-\frac8{9\pi^2}},
$$
if $\varphi_0\neq0$. Since the function $\frac{2-2\cos\varphi}{\varphi^2}$ is non-increasing on $]0,2\pi]$, we would have
$$
\frac{2-2\cos\varphi_0}{\varphi_0^2}<\frac{8}{9\pi^2},
$$
if $\varphi_0$ were greater than $3\pi/2$; thus, we must necessarily have $\varphi_0\in[0,3\pi/2]$. Checking additionally the case $\varphi_0=0$, we obtain that the point $F(\psi_0,\varphi_0)$ belongs to the image of the curve $\tilde F\colon[0,3\pi/2]\to\bR^3$, such that
\begin{align}
\tilde F_1(\varphi)=&\sqrt{\frac{2-2\cos\varphi}{\varphi^2}-\frac8{9\pi^2}};\notag\\
\tilde F_2(\varphi)=&\frac{2\sqrt{2}}{3\pi};\\
\tilde F_3(\varphi)=&2\left(\frac{\varphi-\sin\varphi}{\varphi^2}\right);\notag
\end{align}
where we understand $F_1(0)=\sqrt{1-8/9\pi^2}$ and $F_3(0)=0$.
\end{comment}
%
\beq
\tilde F(\varphi)=\pa{ \sqrt{\frac{2-2\cos\varphi}{\varphi^2}-\frac8{9\pi^2}}, \frac{2\sqrt{2}}{3\pi}, \frac{2\varphi-2\sin\varphi}{\varphi^2}},
\eeq
where we understand $\tilde F(0)=\big(\sqrt{1-8/9\pi^2},\frac{2\sqrt{2}}{3\pi},0\big)$. The image of this curve is the set
\[
\set{(x,2\sqrt{2}/3\pi,z)\in\R^3: x,z\ge0}\cap \der(\B_\rho\cap N_1).
\]
Thus, setting $A^1_\rho=\B_\rho\cap\set{(x,2\sqrt{2}/3\pi,z)\in\R^3:x,z\ge0}$, Gauss-Green's theorem gives
\beq\label{eq:2sect}
\cH^2_{|\cdot|}(A^1_\rho)=2\int_0^{\frac{3\pi}2}\sqrt{\frac{2-2\cos\varphi}{\varphi^2}-\frac8{9\pi^2}}\,\,
\pa{\frac{2\sin\varphi-\varphi\cos\varphi-\varphi}{\varphi^3}}d\varphi,
\eeq
where again $\cH^2_{|\cdot|}(\B_\rho\cap N_1)=4\cH^2_{|\cdot|}(A^1_\rho)$, in view of the symmetries of $\B_\rho$.
Using for instance the computer program {\em Maple}, one can verify that the integral in \eqref{eq:0sect} is strictly less 
than the integral in \eqref{eq:2sect}, hence
\beq\label{eq:strict1}
\gamma_\rho=\frac{\beta^0_\rho}{\beta_\rho}\le \frac{\beta^0_\rho}{\cH^2_{|\cdot|}(\B_\rho\cap N_1)}<1.
\eeq
This conlcudes the proof of Theorem~\ref{th:densless1}.

\subsection{Proof of Theorem~\ref{th:SlessC3}}
Making the same reduction of the previous proof, we can replace $\cS^3_\rho\res\Sigma$ by the perimeter measure $\Per\res\Omega$, for an h-finite perimeter set $E\subset\H$.
The measure theoretic area formula (1.4) of \cite{FSSC8} and Theorem~3.1 of \cite{FSSC5}, joined with Theorem~1.2 of \cite{Mag31}, give the equalities
\[
\Per\res\Omega=\beta^0_\rho\, \cC^3_\rho\res\Sigma=\beta_\rho\,\cS^3_\rho\res\Sigma.
\]
In view of \eqref{eq:strict1}, we have $\cS_\rho^3\res(\Sigma)=\frac{\beta^0_\rho}{\beta_\rho}\;\cC_\rho^3(\Sigma)<\cC_\rho^3(\Sigma)$,
hence the fact that $\cS_\rho^3(T)\le \cC_\rho^3(T)$ on every $T\subset\H$, concludes the proof.

\br\label{re:cCcS3}\rm
When the unit ball with respect to a homogeneous distance $d$ is convex, then 
the equality $\cS^3_d\res\Sigma=\cC^3_d\res\Sigma$ was proved in \cite{FSSC8} for $d=d_\infty$ and general $\G$-rectifiable sets.
The case where $d$ is not required to have a special formula can be recovered following the same arguments of \cite{FSSC8}, 
joined with Theorem~5.2 of \cite{Mag31}.  For the reader's convenience, we will sketch a proof of this equality. 
\er

We follow notation and assumptions in the proof of Theorem~\ref{th:densless1}.
By Theorem~5.2 of \cite{Mag31}, for any horizontal direction $v\in\R^2\times\set{0}$ the following constants 
\beq\label{eq:betad}
\beta(d,v)=\max_{z\in \B_d} \cH_{|\cdot|}^2\lls\B_d\cap zN(v)\rrs\qandq \beta^0(d,v)=\cH^2_{|\cdot|}\lls\B_d\cap N(v)\rrs
\eeq
coincide, where $N(v)=\set{(x,y,z)\in\R^3: \lan(x,y),v\ran=0}$.
By Theorem~3.1 of \cite{FSSC5}, that can be easily extended to any homogeneous distance,
we have 
\[
\Theta^3(\Per\res\Omega,x)=\lim_{r\to0^+}\frac{|\der_HE|(\B(x,r))}{r^3}=\cH^2_{|\cdot|}\lls \B_d\cap N\lls\nu_E(x)\rrs\rrs= \beta^0\lls d,\nu_E(x)\rrs
\]
at $\cS^3_\rho$-a.e.\ $x\in\Sigma$. Thus, by formula (1.4) of \cite{FSSC8} joined with Theorem~1.2 of \cite{Mag31}, we achieve
\[
\Per\res\Omega=\beta^0\lls d,\nu_E(x)\rrs\, \cC^3_0\res\Sigma=\beta\lls d,\nu_E(x)\rrs \cS^3_0\res\Sigma,
\]
where both $\cS_0^3$ and $\cC^3_0$ are constructed with respect to $d$. The equality $\beta(d,\cdot)=\beta^0(d,\cdot)$ finally implies
\beq\label{eq:SC3}
\cS_d^3\res\Sigma=\cC_d^3\res\Sigma.
\eeq

%
%

%
\begin{comment}
In fact, the importance of rectifiability in the Euclidean case can be easily seen. We may consider a purely 1-unrectifiable $A\subset \R^2$ with $\cH^1(A)>0$ and such that $\Theta^{*1}(\cH^1\res A,x)=\frac{1}{2}$ for $\cH^1$ a.e.\ $x\in A$, see 3.3.19 of \cite{Federer69}. 
\end{comment}
%

\vskip2mm
{\bf Acknowledgements.} It is my pleasure to thank Aleksandra Zapadinskaya for her kind contribution in finding the areas of slices of the sub-Riemannian unit ball. I am also grateful to Francesco Serra Cassano for useful comments.

\bibliography{bibtex}{}
\bibliographystyle{plain}

\end{document}